\documentclass[onethmnum,oneeqnum]{siamltex}

\usepackage{amsfonts}
\usepackage{amsmath}

\newcommand{\dd}{\, {\mathrm d}}
\newcommand{\QG}{Q^{\mathrm G}}

\newcommand{\QCC}{Q^{\mathrm{CC}}}

\newcommand{\Var}{\mathop{\mathrm{Var}}}

\title{Error Bounds for the Numerical Integration\\ of Functions with Limited Smoothness}
\author{Kai Diethelm\footnotemark[2]\ \footnotemark[3]}

\begin{document}
\maketitle

\renewcommand{\thefootnote}{\fnsymbol{footnote}}
\footnotetext[2]{GNS Gesellschaft f\"ur numerische Simulation mbH, Am Gau\ss berg 2, 38114 Braunschweig, Germany, \tt diethelm@gns-mbh.com}
\footnotetext[3]{Institut Computational Mathematics, Technische Universit\"at Braunschweig, Fallersleber-Tor-Wall 23, 38100 Braunschweig, Germany, \tt k.diethelm@tu-bs.de}

\begin{abstract}
  Recently, Trefethen (SIAM Review 50 (2008), 67--87) and 
  Xiang and Bornemann (SIAM J. Numer.\ Anal.\ 50 (2012), 2581--2587)
  investigated error bounds for $n$-point Gauss and Clenshaw-Curtis
  quadrature for the Legendre weight with integrands having
  limited smoothness properties. Putting their results
  into the context of classical quadrature theory, we
  find that the observed behaviour is by no means surprising and that
  it can essentially be proved for a very large class of quadrature
  formulas with respect to a broad set of weight functions.
\end{abstract}

\begin{keywords}
  Numerical integration; Gauss quadrature; Clenshaw-Curtis quadrature; error bound; 
  bounded variation; weight function.
\end{keywords}

\begin{AM}
  41A55
\end{AM}

\pagestyle{myheadings}
\thispagestyle{plain}
\markboth{K. DIETHELM}{INTEGRATION OF FUNCTIONS WITH LIMITED SMOOTHNESS}

%% \section{Introduction}
The numerical approximation of the weighted integral
\begin{equation}
  \label{eq:def-i}
  I_w[f] := \int_{-1}^1 w(x) f(x) \dd x
\end{equation}
% by means of interpolatory quadrature formulas,
with a nonnegative integrable weight function $w$ and $f$
possessing certain smoothness properties is a classical topic of
research \cite{Br,BP} that recently attracted new attention. 
Specifically, Trefethen \cite{Tr08} looked at the $n$-point formulas
of Gauss type $\QG_n$ and Clenshaw-Curtis type $\QCC_n$, respectively,
for $w \equiv 1$ and showed assuming that $f^{(s-1)}$ is
absolutely continuous for some $s \in \{1,2,3,\ldots\}$, $\Var f^{(s)} <
\infty$ (where $\Var g$ denotes the total
variation of $g$ over the interval $[-1,1]$) and that $f^{(s)}$ is sufficiently
well behaved near the points $\pm1$, 
cf.\  \cite[p.\ 75 and Theorems 4.5 and 5.1]{Tr08}, that 
\begin{equation}
  \label{eq:tr}
    I_w[f] - \QCC_n[f] = O(n^{-s}) 
    \quad \mbox{ and } \quad
    I_w[f] - \QG_n[f]  = O(n^{-s}).
\end{equation}
Later, without requiring Trefethen's
assumptions about the behaviour of the integrand near the boundary,
Xiang and Bornemann \cite{XB} improved the bounds to  
\begin{equation}
  \label{eq:xb}
    I_w[f] - \QCC_n[f] = O(n^{-s-1}) 
    \quad \mbox{ and } \quad
    I_w[f] - \QG_n[f]  = O(n^{-s-1});
\end{equation}
however, their proof for the Gauss formula needs the additional
restriction that $s \ge 2$.

We now want to place these results into
the context of classical quadrature theory. 
In this connection, a \emph{quadrature formula} is a linear functional 
$Q$ of the form $Q[f] := \sum_{j=1}^n a_j f(x_j)$
with \emph{nodes} $-1 \le x_1 < x_2 < \ldots < x_n \le 1$ and
\emph{weights} $a_j \in \mathbb R$. Following the common
terminology in numerical integration \cite{BP}, we say that such
a quadrature formula is \emph{positive} if $a_j \ge 0$ for all $j$,
and a formula is \emph{interpolatory} if $Q[p] =
I_w[p]$ for all $p \in \mathcal P_{n-1}$, the set of all polynomials
of degree not exceeding $n-1$.

%% \section{Classical Results}
To achieve our goal we define the function classes 
\[
  V_0 :=  \{f \! : \! [-1,1] \to \mathbb R : \Var f < \infty \}
\]
and
\[
  V_s :=  \{f \! : \! [-1,1] \to \mathbb R : f^{(s-1)} \mbox{ is absolutely continuous and }
                                        \Var f^{(s)} < \infty \}
   \mbox{ for } s \in \mathbb N 
\]
and recall an
important result of Freud \cite[Satz II]{Fr55}
(see also \cite[\S III.4]{Fr} for the case $s=0$)
who proved, using an
argument based on one-sided polynomial approximation in the $L^1$ norm, 
the following statement: 

\begin{theorem}
  Let $Q_n$ be a positive interpolatory quadrature formula with
  $n$ nodes for the weight $w$, and assume $w(x) \le M (1-x^2)^{-1/2}$ for some $M > 0$. 
  Then,
  $$ |I_w[f] - Q_n[f]| \le 5 M ((s+2)\pi)^{s+1} (s!)^{-1} n^{-s-1} \Var f^{(s)} $$
  whenever $f \in V_s$ for some $s \in \mathbb N_0$.
\end{theorem}

Freud's result has some immediate consequences: 
\begin{remunerate}
\item The sharper bounds (\ref{eq:xb}) hold for all functions $f$
  with an absolutely continuous derivative of order $s-1$ and 
  $\Var f^{(s)} < \infty$ for all $s \in \mathbb N_0$ 
  (again without additional assumptions near
  the boundary of the interval).
\item Bounds of the same order also hold for a very large class of weight
  functions and many other quadrature formulas including, e.g., 
  the Gauss and Radau formulas for all admissible weight functions \cite{Br,BP},
  the formulas of Clenshaw-Curtis \cite{CC} and their close relatives
  due to Filippi \cite{Fi} and Polya \cite{Po}, at least for  
  the standard weight function $w \equiv 1$, and the Gauss-Kronrod
  formulas \cite{Ga,Kr} for the ultraspherical weight functions 
  $w_\lambda(x) = (1-x^2)^{\lambda - 1/2}$ where $\lambda \in
  [0,1] \cup \{3\}$ (cf.\ also \cite{Pe99}).
\end{remunerate}

Other traditional techniques from the theory of numerical integration
provide additional insight. For example, 
positive interpolatory quadrature
formulas are not the only numerical integration methods admitting an
error estimate of the form (\ref{eq:xb}) for\linebreak functions $f \in V_s$: 
At least for the Legendre weight function
$w \equiv 1$, compound quadrature schemes exhibit the same behaviour. 
To see this, let $Q$ be an arbitrary quadrature formula given by 
$Q[f] = \sum_{j=1}^\ell a_j f(x_j)$. We
then subdivide the basic interval $[-1,1]$ into $n$ subintervals
of length $2/n$, affinely transform $Q$ to each of these subintervals
and add up the resulting formulas, thus obtaining the $n$-fold
\emph{compound quadrature formula} with respect to the
\emph{elementary formula} $Q$, denoted and defined by
$$
  Q^{(n)}[f] := \frac1n \sum_{\nu=1}^n \sum_{j=1}^\ell a_j f\left(-1 + \frac1n (x_j + 2 \nu - 1) \right).
$$
Using the Peano kernel representation of the associated error
functional (cf.\ \cite[Satz 17, Satz 93 and Satz 97]{Br} or \cite[\S\S4.2 and 4.4]{BP}), 
one can immediately conclude

\begin{theorem}
  \label{thm:compound}
  Let $w \equiv 1$, let $Q^{(n)}$ be the $n$-fold compound quadrature
  formula for the elementary formula $Q$, and assume that $Q[p] =
  I_w[p]$ for all $p \in \mathcal P_m$ with some $m \in \mathbb N_0$.
  Then, for every $s \in \{0,1,\ldots,m\}$ there exists a constant $C$
  (depending only on $Q$ and $s$) such that,  for all $f \in V_s$,
  $$ 
    |I_w[f] - Q^{(n)}[f]| \le C n^{-s-1} \Var f^{(s)}. 
  $$
\end{theorem}
\vskip-1em

Thus, Xiang's and Bornemann's improvements of Trefethen's estimates
for the error of Gaussian or Clenshaw-Curtis quadrature for integrands
having an $s$th derivative of bounded variation fit into a more
general picture. 

%% \section{Conclusion}
%% 
%% Using tools from classical quadrature theory, we have seen that
%% Xiang's and Bornemann's improvements of Trefethen's estimates
%% for the error of Gaussian or Clenshaw-Curtis quadrature for integrands
%% having an $s$th derivative of bounded variation fits into a more
%% general picture. Specifically, error estimates of the same order of
%% magnitude can be seen to hold (for all $s$) for very many other positive
%% interpolatory quadrature methods and a large class of weight
%% functions. Moreover, for the Legendre weight function, these error
%% estimates hold for restricted values of $s$ also for compound
%% quadrature methods.

\end{document}